\documentclass[final]{article}
\usepackage{amscd}
\usepackage{amsxtra}
\usepackage{amssymb, amsmath}
\usepackage{amsthm}

\theoremstyle{plain}
\newtheorem{theorem}{Theorem}[subsection]
\newtheorem{proposition}[theorem]{Proposition}

\newtheorem{corollary}[theorem]{Corollary}
\newtheorem{definition}[theorem]{Definition}
\newtheorem{remark}[theorem]{Remark}

\begin{document}

\newcommand{\U}{\ensuremath{\mathfrak{u}_{1}}}
\newcommand{\ad}{\ensuremath { ad(\mathfrak{u}_{1})}}
\newcommand{\Oo}{\ensuremath{\varOmega^{0}(ad(\mathfrak{u}_{1}))}}
\newcommand{\Ok}{\ensuremath {\varOmega^{1}(ad(\mathfrak{u}_{1}))}}

\newcommand{\s}{\ensuremath{\mathcal{S}}}
\newcommand{\csa }{\ensuremath {\mathcal{S_{\alpha}}}}
\newcommand{\cs}{\ensuremath{\mathcal{S^{c}}}}
\newcommand{\vsa}{\ensuremath{\varGamma(\mathcal{S^{+}_{\alpha}})}}
\newcommand{\la }{\ensuremath{{\mathcal{L}}_{\alpha}}}
\newcommand{\sla }{\ensuremath{{\mathcal{L}}^{1/2}_{\alpha}}}

\newcommand{\ca}{\ensuremath{\mathcal{C}_{\alpha}}}
\newcommand{\Aa}{\ensuremath {\mathcal{A}_{\alpha}}}
\newcommand{\Q}{\ensuremath{\mathcal{A}_{\alpha}\times_{\mathcal{G}_{\alpha}}
\varGamma (S^{+}_{\alpha})}}
\newcommand{\wQ}{\ensuremath{\mathcal{A}_{\alpha}\times_{\widehat{\mathcal{G}}_{\alpha}}
\varGamma (S^{+}_{\alpha})}}
\newcommand{\B}{\ensuremath{\mathcal{B}_{\alpha}}}
\newcommand{\wB}{\ensuremath{\widehat{\mathcal{B}}_{\alpha}}}

\newcommand{\G}{\ensuremath{\mathcal{G}_{\alpha}}}
\newcommand{\wG}{\ensuremath{\widehat{\mathcal{G}}_{\alpha}}}

\newcommand{\lda}{\lambda}

\newcommand{\spinc}{\ensuremath{Spin^{c}_{4}}}
\newcommand{\spin }{\ensuremath{Spin_{4}}}

\newcommand{\Z }{\ensuremath{\mathbb {Z}}}
\newcommand{\R }{\ensuremath{\mathbb {R}}}
\newcommand{\C }{\ensuremath {\mathbb {C}}}
\newcommand{\N }{\ensuremath{\mathbb {N}}}
\newcommand{\iso }{\ensuremath {\thickapprox }}

\newcommand{\sw}{\ensuremath {\mathcal{S}\mathcal{W}_{\alpha}}}
\newcommand{\y}{\ensuremath {\mathcal{Y}\mathcal{M}}}
\newcommand{\yp}{\ensuremath {\mathcal{Y}\mathcal{M}^{+}}}
\newcommand{\yn}{\ensuremath {\mathcal{Y}\mathcal{M}^{-}}}
\newcommand{\Cf }{\ensuremath {\mathcal{C}_{\alpha}}}
\newcommand{\w }{\ensuremath {\omega}}
\newcommand{\cx}{\ensuremath {C_{X}}}

\title{On the Existence of Critical Points to the Seiberg-Witten functional}
\author{ Celso M. Doria \\ UFSC - Depto. de Matem\'atica }
\maketitle
\begin{abstract}
Let X be a closed smooth 4-manifold. In the Theory of the Seiberg-Witten 
Equations\footnote{ MSC 58J05 , 58E50}\footnote{connections,gauge fields,
4-manifolds }, the configuration space is $\mathcal{C}_{\alpha}=
\mathcal{A}_{\alpha}\times \varGamma (S^{+}_{\alpha})$,
where $\mathcal{A}_{\alpha}$ is a space of $\mathfrak{u}_{1}$-connections
defined on a complex line bundle over X and $\varGamma(S^{+}_{\alpha})$ is the 
space of sections of the positive complex spinor bundle over X. The original 
$\sw$-equations are $1^{st}$-order PDE fitting into a variational principle 
$\sw:\Aa \times\csa\rightarrow \R$, which is invariant by the group action 
of (Gauge Group) $\mathcal{G}_{\alpha}=Map(X,U_{1})$ and 
satisfies the Palais-Smale Condition, up to gauge equivalence. The Euler-Lagrange equations of
the functional $\sw$ are $2^{nd}$-order PDE and the solutions of the original
 $\sw$-equations are  stable critical points. Our aim is to prove the existence of 
solutions to the Euler-Lagrange equations of the functional $\sw$ by the method of the 
Minimax Principle.

\end{abstract}

\section{\bf{Introduction}}

\indent Although the physical meaning of the Seiberg-Witten equations ($\sw$-eq.)
 is yet to be discovered, the mathematical meaning is rather deep and highly
 efficient to understand one of the most basic phenomenon of differential topology in four
dimension, namely, the existence of  non-equivalent differential smooth
structures on the same underlying topological manifold.
The Seiberg-Witten equations arised through the  ideas of duality described in
Witten ~\cite{Wi94}. It is conjectured that the Seiberg-Witten equations
are dual to Yang-Mills equations ($\y$-eq.). The duality is at the  quantum
 level, since one of its necessary condition is the equality of the
 expectations values for the dual theories. In topology, this means that fixed a
 4-manifold its Seiberg-Witten invariants are equal to Donaldson invariants. A
 good reference for $\sw$-eq. is ~\cite{Do96} 

 Originally, the $\sw$-equations were 1$^{st}$-order differential equations
without a variational principle associated to them. In ~\cite{Wi94}, Witten
used some special identities to obtain a integral useful to prove that the
moduli spaces of the theory were empty, but a finite number of
them. Jost-Peng-Wang, in
~\cite{JPW96}, used this integral to define a functional, which we refer
as the $\sw$ functional; their main result is to prove that the functional
satisfies the Palais-Smale (PS)-condition up to a gauge equivalence.  

In our context, we refer to the $\sw$-equation as the
Euler-Lagrange equation of the $\sw$ functional. These equations
are 2$^{nd}$-order differential equations and from now on they will be called 
$\sw$-equations. 

 Since the  (PS)-condition is satisfied in the quotient space, the main aim is 
to  describe the weak homotopy type of the moduli space $\Q$ in order to prove 
the existence of solutions to the $\sw$-equation.

\section{\bf{Basic Set Up }}

From a duality principle appliable to SUSY theories in Quantum Field
Theory, Seiberg-Witten, in ~\cite{SW94}, discovered a nice coupling of the 
self-dual(\emph{SD}) equation, of a $U_{1}$
 Yang-Mills  Theory, to the Dirac$^{c}$ equation. In order to describe this
 coupling it is necessary  a particular isomorphism relating the space
 $\varOmega^{2}_{+}(X)$, of \emph{self-dual} 2-forms, and the bundle 
 $End^{0}(\csa^{+}$) (~\cite{JM96}).

  The space of $Spin^{c}$-strutures on X is identified as
 $$Spin^{c}(X)=\{\alpha\in H^{2}(X,\mathbb{Z})\mid w_{2}(X) = \alpha \ mod
 \ 2\}.$$ 

 For each $\alpha \in Spin^{c}(X)$, there is a representation 
$\rho_{\alpha}:SO_{4}\rightarrow \C l_{4}$, induced by a \emph{$Spin^{c}$} representation,
 and consequently, a pair of vector bundles 
$(\csa^{+},\la)$ over X (see ~\cite{LM89}), where
\begin{itemize}
\item $\csa = P_{SO_{4}}\times_{\rho_{\alpha}} V = \csa^{+}\oplus
  \csa^{-}$.\\
 The bundle  $\csa^{+}$ is the positive complex spinors bundle (fibers are 
$Spin^{c}_{4}-modules$ isomorphic to $\C^{2}$) 
\item $\la=P_{SO_{4}}\times_{det(\alpha)} \C$.\\
 It is called the \emph{determinant line bundle} associated to the
 $Spin^{c}$-struture $\alpha$. ($c_{1}(\la)=\alpha$)
\end{itemize}

Thus, given $\alpha \in Spin^{c}(X)$ we associate a pair of bundles

$$\alpha \in Spin^{c}(X) \quad \rightsquigarrow \qquad(\la,\csa^{+})$$

From now on, we  fixed 
\begin{itemize}
\item a Riemannian metric g over X 
\item a Hermitian structure h on $\csa$.
\end{itemize}

\begin{remark} Let $E\rightarrow X$ be a vector bundle over X (~\cite{EL80});
\begin{enumerate}
\item The space of sections of E (usually denoted by $\varGamma(E)$) is denoted by $\varOmega^{0}(E)$;
\item The space of p-forms ($1\le p\le 4$) with values in E is denoted by $\varOmega^{p}(E)$.
\item For each fixed covariant derivative $\triangledown$
 \footnote{on E, connection 1-form A $\leftrightarrow$ $\triangledown^{A}$ covariant derivative}
 on E, there is a 1$^{st}$-order differential operator 
$d^{\triangledown}:\varOmega^{p}(E)\rightarrow \varOmega^{p+1}(E)$  
\end{enumerate}
\end{remark}

For each class $\alpha\in Spin^{c}(X)$ corresponds  a $U_{1}$-principal bundle over X, denoted
  $P_{\alpha}$, with $c_{1}(P_{\alpha})=\alpha$. Also, we consider the adjoint bundles 

$$Ad(U_{1})=P_{U_{1}}\times_{Ad}U_{1}\quad ad(\U)=P_{U_{1}}\times_{ad}\U.$$

 $Ad(U_{1})$ is  a fiber bundle with fiber $U_{1}$, and $ad(\mathfrak{u}_{1})$ is 
a vector bundle with fiber isomorphic to the Lie Algebra $\mathfrak{u}_{1}$. 

Once a covariant derivative is fixed on $\ad$, it induces the sequence

\[
 \begin{CD}
  \Oo @>d^{\triangledown}>> \varOmega^{1}(ad(\mathfrak{u}_{1}))
  @>d^{\triangledown}>> \varOmega^{2}(ad(\mathfrak{u}_{1})) @>d^{\triangledown}>> (*)
\end{CD} 
\]

\[
\begin{CD}
(*)  @>d^{\triangledown}>> \varOmega^{3}(ad(\mathfrak{u}_{1})) @>d^{\triangledown}>>   \varOmega^{4}(ad(\mathfrak{u}_{1})) 
 \end{CD} 
  \]

\vspace{05pt}

   The 2-form of curvature $F_{\triangledown}$, induced by the connection $\triangledown$,
 is the operator
 
$$F_{\triangledown} = d^{\triangledown}\circ d^{\triangledown}:\Oo \rightarrow \varOmega^{2}(ad(\mathfrak{u}_{1}))$$

\vspace{05pt} 
  
Since $Ad(U_{1})\sim X\times U_{1}$ and $\ad \sim X\times \U$, the spaces  $\Oo$ and 
$\varGamma(Ad(U_{1}))$ are identified, respectively, to the spaces
$\varOmega^{0}(X,i\R)$ and $Map(X,U_{1})$. It is well known from the theory (see in ~\cite{DK91}) that a
$\mathfrak{u}_{1}$-connection defined on $\la$ can be identified with
 a section of the vector bundle $\Ok$, and a Gauge transformation with a section of the 
bundle $Ad(U_{1})$. 

Given a vector bundle E over (X,g), endowed with a metric and a covariant derivative 
$\triangledown$, we define the Sobolev Norm of a section $\phi \in \varOmega^{0}(E)$ as 

$$\mid\mid\phi\mid\mid_{L^{k,p}}=\sum_{\mid i \mid =0}^{k}
(\int_{X}\mid\triangledown^{i}\phi\mid^{p})^{\frac{1}{p}}$$

and the Sobolev Spaces of sections of E as
$$L^{k,p}(E)=\{\phi\in \varOmega^{0}(X,E) \mid \quad \mid\mid\phi\mid\mid_{L^{k,p}} < \infty\}$$

Now, consider the spaces
\begin{itemize}
\item $\mathcal{A}_{\alpha}= L^{1,2}(\Oo)$ 
\item $\varGamma(\csa^{+})$ = $L^{1,2}(\varOmega^{0}(X,\csa^{+})$  
\item $\ca = \mathcal{A}_{\alpha}\times \varGamma(\csa^{+})$
\item $\G= L^{2,2}(X,U_{1})=L^{2,2}(Map(X,U_{1}))$
\end{itemize}
 
The space $\G$ is the Gauge Group acting on $\ca$ by the action

\begin{equation}
\G\times\ca \rightarrow \ca ;\quad
 (g,(A,\phi ))\rightarrow (g^{-1}dg + A, g^{-1}\phi)
\end{equation}

\vspace{05pt}

Since we are in dimension 4, the vector bundle $\varOmega^{2}(ad(\mathfrak{u}_{1}))$ admits a decomposition

\begin{equation}
\varOmega^{2}_{+}(ad(\mathfrak{u}_{1}))\oplus \varOmega^{2}_{-}(ad(\mathfrak{u}_{1}))
\end{equation}

\noindent in  seld-dual (+) and anti-self-dual (-) parts (~\cite{DK91}).

\vspace{05pt}

 The 1$^{st}$-order (original) \emph{Seiberg-Witten} equations are defined  over the configuration space
$\ca =\mathcal{A}_{\alpha}\times \varGamma(\csa^{+})$ as

\begin{equation}\label{E:03}
\begin{cases}
D^{+}_{A}(\phi )= 0,\\ 
F^{+}_{A} = \sigma (\phi)
\end{cases}
\end{equation}  
where
\begin{itemize}
\item $D^{+}_{A}$ is the $Spinc^{c}$-Dirac operator defined on $\vsa$;
\item Given $\phi \in \varGamma(\csa^{+})$, the quadratic form 

\begin{equation}\label{E:QF}
\sigma (\phi)=\phi\otimes\phi^{*} - \frac{\mid\phi\mid^{2}}{2}.I
\qquad\Rightarrow \qquad 
\sigma (\phi)\in End^{0}(\csa^{+})
\end{equation}

\vspace{05pt}

\noindent performs the coupling of the \emph{ASD}-equation with
the $Dirac^{c}$ operator. \\
Locally, if $\phi = (\phi_{1},\phi_{2})$, then the quadratic form
$\sigma(\phi)$ is written as 

\begin{equation*}
 \sigma(\phi) = \left(
 \begin{matrix}
 \frac{\mid\phi_{1}\mid^{2}-\mid\phi_{2}\mid^{2}}{2} & \phi_{1}.\Bar{\phi_{2}} \\
 \phi_{2}.\Bar{\phi_{1}} &\frac{ \mid\phi_{2}\mid^{2} - \mid\phi_{1}\mid^{2}}{2}
 \end{matrix}
 \right)
 \end{equation*}
\end{itemize}

The set of solutions of equations (~\ref{E:03}) can be described as the inverse image 
$\mathcal{F}^{-1}(0)$ by the map  
$\mathcal{F}_{\alpha} :\ca \rightarrow \varOmega^{2}_{+}(X)\oplus
 \varGamma(\csa^{-})$, defined as

$$\mathcal{F}_{\alpha}(A,\phi) = (F^{+}_{A} - \sigma (\phi), D^{+}_{A}(\phi) )$$ 

\vspace{05pt}
 
The $\sw$-equations are $\G$-invariant.

\section{\bf{A Variational Principle for the Seiberg-Witten Equation}}

The lack of a natural \emph{Lagrangean} let us to consider, in a pure formally
way, the functional

\begin{equation}\label{E:SW00}
SW (A,\phi) =\frac{1}{2} \int_{X}\{\mid F^{+}_{A} - \sigma(\phi)\mid^{2} + \mid
D^{+}_{A}(\phi)\mid^{2}\}dv_{g}
\end{equation}

\noindent The next identities, which proofs are standard, are useful to 
expand the functional (~\ref{E:SW00})

\begin{proposition}
For each $\alpha \in Spin^{c}(X)$, let $\mathcal{L}_{\alpha}$ be the
determinant line  bundle associated to $\alpha$ and (A,$\phi$) $\in\ca$. Also, assume that 
$k_{g}$=scalar curvature of (X,g). Then,

\begin{enumerate}
\item $<F^{+}_{A},\sigma(\phi)> = \frac{1}{2} <F^{+}_{A}.\phi ,\phi >$
\vspace{04pt}
\item $<\sigma(\phi ),\sigma (\phi )> =\frac{1}{4}\mid \phi \mid^{4} $
\vspace{04pt}
\item Weitzenb$\ddot{o}$ck formula $$ D^{2}\phi = \triangledown^{*}\triangledown \phi +
  \frac{k_{g}}{4}\phi +\frac{F_{A}}{2}.\phi $$

\item $\sigma(\phi)\phi = \frac{\mid\phi\mid^{2}}{2}\phi$
\vspace{04pt}
\item $c_{2}(\la\oplus \la)=\int_{X}F_{A}\wedge F_{A} $
\vspace{04pt}
\item $\mid F^{+}_{A}\mid^{2} = \frac{1}{2}\mid F_{A}\mid^{2} - 4\pi^{2}\alpha^{2}$
\end{enumerate}
\end{proposition}

Consequently, after expanding the functional (~\ref{E:SW00}), we  get the expression

\begin{equation}\label{E:SW01}
SW(A,\phi) = \int_{X}\{\frac{1}{4}\mid F_{A}\mid^{2} + \mid \triangledown^{A} \phi
\mid^{2} + \frac{1}{8}\mid \phi\mid^{4} + \frac{1}{4}<k_{g}\phi , \phi> \}dv_{g}- 2\pi^{2}\alpha^{2}
\end{equation}

\vspace{05pt} 

\begin{definition}
For each $\alpha\in Spin^{c}(X)$, the Seiberg-Witten Functional is the functional 
 $\sw:\ca\rightarrow\R$ given by
 \begin{equation}\label{E:SW02}
\sw(A,\phi) = \int_{X}\{\frac{1}{4}\mid F_{A}\mid^{2} + \mid \triangledown^{A} \phi
\mid^{2} + \frac{1}{8}\mid \phi\mid^{4} + \frac{1}{4}<k_{g}\phi , \phi> \}dv_{g}
\end{equation}
where $k_{g}$= scalar curvature of (X,g).
\end{definition}

Let $k_{g,X}=min_{x\in X}k_{g}$ and 

\begin{equation}\label{E:SC}
k^{-}_{g,X}=min\{0, - k^{\frac{1}{2}}_{g,X}\}
\end{equation}

\begin{remark}
  \begin{enumerate}
\item Since X is compact and $\mid\mid \phi \mid\mid_{L^{4}}<\mid\mid \phi
    \mid\mid_{L^{1,2}}$, the functional is well defined on $\ca$ ,
\item The $\sw$-functional \emph{(\ref{E:SW02})} is Gauge invariant.
\item From (~\ref{E:SW00}) and (~\ref{E:SW01}), it follows that 

$$\sw (A,\phi)-2\pi^{2}\alpha^{2}\ge 0$$

\noindent Thus, the $\sw$-functional is bounded below by $2\pi^{2}\alpha^{2}$, where 
$$\alpha^{2}= Q_{X}(\alpha , \alpha)$$
 ($Q_{X}:H^{2}(X,\Z)\times H^{2}(X,\Z)\rightarrow \Z$ is the
 intersection form of X). 
\vspace{04pt}
\item The set of classes in $Spin^{c}(X)$, such that there exits  $(A,\phi)\in \ca$ attaining the 
minimum value $2\pi^{2}\alpha^{2}$, is finite (~\cite{JM96}). More precisely, the minimum value
 $2\pi^{2}\alpha^{2}$ is attained for those $\alpha \in Spin^{c}(X)$ such that

$$\alpha^{2}\le \frac{1}{4\pi^{2}}.(\frac{k^{-}_{g,X}}{8}.vol(X) + 2\chi(X) +3\sigma_{X})$$

where $\chi(X)=\text{Euler characteristic class of X}$ and $\sigma_{X}$=signature of $Q_{X}$.

  \end{enumerate}
\end{remark}

\begin{proposition}
The Euler-Lagrange equations of the $\sw$-functional \emph{(\ref{E:SW02})} are

\begin{equation}\label{E:09}
\Delta_{A} \phi + \frac{\mid \phi \mid^{2}}{4}\phi +\frac{k_{g}}{4}\phi = 0 
\end{equation}
\begin{equation}\label{E:10}
d^{*}F_{A} + 4\Phi^{*}(\triangledown^{A}\phi)=0
\end{equation}
where $\Phi:\varOmega^{1}(\mathfrak{u}_{1})\rightarrow\varOmega^{1}(\csa^{+})$
\end{proposition}

\begin{proof}
Let $\phi_{t}:(-\epsilon,\epsilon)\rightarrow \varGamma(\csa^{+})$ be a small
perturbation of a critical point $(A,\phi)\in \ca$, such that $\phi_{0}=\phi$
and $\frac{d\phi_{t}}{dt}\mid_{t=0}=\Lambda\in \vsa$. The derivative of the function \\ 
$\sw(t)=\sw (A,\phi_{t})$, at
$(A,\phi)$, is given by

$$\frac{d\sw(t)}{dt}\mid_{t=0}=d(\sw)_{(A,\phi)}.\Lambda = \int_{X}Re\{<\Delta_{A} \phi + \frac{\mid
\phi \mid^{2}}{4}\phi + \frac{k_{g}}{4}\phi, \Lambda>\}dx ,$$

\noindent where $\Delta_{A} = (\triangledown^{A})^{*}.\triangledown^{A}$.

\noindent If for all $\Lambda \in \varGamma(S^{+})$, we have that
$d(\sw)_{(A,\phi)}.\Lambda =0$

\noindent then

  $$\Delta_{A} \phi + \frac{\mid \phi \mid^{2}}{4}\phi +\frac{k_{g}}{4}\phi = 0 $$

\noindent The second equation is obtained by considering a smooth curve \\
$A_{t}:(-\epsilon,\epsilon)\rightarrow \mathcal{A}_{\alpha}$ given by 
$A_{t}=A + t\Theta$, where $\Theta\in \varOmega^{1}(\mathfrak{u}_{1})$. 
It follows from the first order approximations 
$$F_{A+t\Theta}=F_{A} + t d_{A}\Theta +
o(t^{2})$$

 \noindent and 

$$\triangledown^{A + t\Theta}\phi = \triangledown^{A}\phi + t\Theta(\phi)+ o(t^{2}),$$

\noindent  that

$$d\sw .\Theta = \frac{1}{4}\int_{X}\{<F_{A},d_{A}\Theta> + 4<\triangledown^{A}(\phi),
\Phi(\Theta)>,$$ 

\noindent where $\Phi:\varOmega^{1}(\mathfrak{u}_{1})\rightarrow
\varOmega^{1}(\csa^{+})$ is the linear operator $\Phi(\Theta) = \Theta(\phi)$

\noindent If $(\phi , A)$ is a critical point of $\sw$, then for all 
$\Theta \in \varGamma^{1}(\mathfrak{u}_{1})$

$$d\sw .\Theta = \frac{1}{4}\int_{X}<d^{*}_{A}F_{A} +
4\Phi^{*}(\triangledown^{A}\phi),\Theta> = 0$$

\noindent Hence, 

$$d^{*}_{A}F_{A} + 4\Phi^{*}(\triangledown^{A}\phi)=0$$

\end{proof}

\begin{remark}
Locally, in a orthonormal basis $\{\eta^{i}\}_{1\le i\le 4}$ of $T^{*}X$, the 
operator $\Phi^{*}$ can be written as 

$$\Phi^{*}(\triangledown^{A}\phi)=
\sum_{i=1}^{4}<\triangledown^{A}_{i}\phi,\phi>\eta^{i}, \quad\text{where}\quad 
\triangledown^{A}_{i}=\triangledown^{A}_{X_{i}}\quad (\eta_{i}(X_{j})=\delta_{ij})$$

\end{remark}

The regularity of the solutions of (~\ref{E:09}) and (~\ref{E:10}) was studied by Jost-Peng-Wang 
 in (~\cite{JPW96}). They observed that the $L^{\infty}$ estimate of a solution $\phi$, 
already known to be satisfied by the stable critical points, is also obeyed by the non-stable critical 
points of $\sw$. The estimate is the following; 
 
\begin{proposition}
If $(A,\phi)\in \ca$ is a solution of (\ref{E:09}) and (\ref{E:10}), then

$$\mid\mid \phi \mid\mid_{\infty}\le k^{-}_{g,X}$$

where $k^{-}_{g,X}=max_{x\in X}\{0, - k^{\frac{1}{2}}_{g}(X)\}$ 
\end{proposition}

As a consequence of the estimate above, if the Riemannian metric g on X has  non-negative 
scalar curvature then the only solutions are $(A,0)$, where

$$d^{*}F_{A}=0$$

Since X is compact,

$$d^{*}F_{A}=0\quad\Leftrightarrow\quad\Delta_{A}F_{A}=0\quad (F_{A}\quad\text{is harmonic})$$
 
It follows from  the formula to the 1$^{st}$-Chern class

$$\alpha([\Sigma])=\frac{1}{2\pi i}\int_{[\Sigma]}F_{A},$$

\noindent for all class $[\Sigma]\in H^{2}(X,\R)$,  that $F_{A}$ is the
only harmonic representative of the De Rham class of  $\alpha$. 

 If $(A,0)$ is a solution of the 1$^{st}$-order $\sw$-equation (minimum for $\sw$), then 
 $F^{+}_{A}=0$ and $\sw(A,0)=2\pi^{2}\alpha^{2}$. It is known (~\cite{DK91}) that if 
$b^{+}_{2}>1$, then such solutions do not exists for a dense set of the space of metrics 
on X. Therefore, whenever  $b^{+}_{2}>1$, there is a dense set of metrics on X such 
that 

$$\sw(A,0)>2\pi^{2}\alpha^{2}.$$

Although, for each $\alpha\in Spin^{c}(X)$, the functional always attains its minimum value 
(~\ref{Th:M}), it may happens that  

$$\inf_{(A,\phi)\in\ca}\sw(A,\phi)>2\pi^{2}\alpha^{2},$$

\noindent since 

$$\inf_{(A,\phi)\in\ca}\sw(A,\phi) = 2\pi^{2}\alpha^{2}$$

\noindent just for a finite subset of $Spin^{c}(X)$.
 
 In the Euclidean $\R^{4}$, the only solution to the equations (~\ref{E:09}) 
and (~\ref{E:10}), up to gauge equivalence, is the trivial one (0,0). 
 
In ~\cite{JPW96}, Jost-Peng-Wang studied the analytical properties of the
$\sw$-functional. They proved that the Palais-Smale Condition, up to gauge equivalence,
 is satisfied. Whenever the quotient space is a smooth manifold, one may  use 
\emph{Minimax Principle}  to prove the existence of solutions to 
(~\ref{E:09}) and (~\ref{E:10}). (see ~\cite{Pa68}).

Since the $\sw$-funtional is $\G$-invariant,  it induces a functional
on the space $\Q$; the quotient space by the Gauge Group action. However, the space
$\Q$ isn't a manifold since the action isn't free.

The $\G$-action on $\ca$ has non-trivial isotropic groups once all elements in $\ca$ 
are fixed by the action of the constant maps $g:X\rightarrow U_{1}$. In fact, 
\begin{enumerate}
\item If $(A,\phi)\ne (0,0)$ $\Rightarrow$ $G_{(A,\phi)}\overset{\text{iso}}{\simeq} U_{1}$, since

$$g.(A,\phi)=A\Leftrightarrow g^{-1}dg=0\Leftrightarrow g\quad\text{is constant}$$

\item $G_{(0,0)}\overset{\text{iso}}{\simeq} \G$
\end{enumerate}  

Therefore, we consider the Gauge Group

$$\wG =\frac{\G}{\{g:X\rightarrow U_{1}\mid\text{ g=constant}\}} \overset{\text{iso}}{\simeq}
\G\diagup U_{1}$$

\vspace{05pt}

If $\alpha\ne 0$, then the $\wG$-action on $\ca$ is free.
Ignoring the case of the trivial bundle ($\alpha =0$); from now on, instead of the $\G$-action, 
we consider on $\ca$ the  $\widehat{\G}$-action, and so, the quotient space $\wB$ is a manifold.
In the case of the trivial bundle we get a orbifold.

Let's consider $\sw:\wQ\rightarrow\R$ as the induced functional. The Palais-Smale Condition proved
by Jost-Peng-Wang can be written in the following way;

\begin{proposition}
(~\cite{JPW96}) Consider a sequence $[(A_{n},\phi_{n})]\in \Q$ satisfying the conditions
\begin{enumerate}
\item $d(\sw)_{[(A_{n},\phi_{n})]}\rightarrow 0$ strongly in $\wQ$,
\item $\sw(A_{n},\phi_{n})\le c$ for $n\in\N$. 
\end{enumerate}
So, there exists a subsequence $[(A_{n_{k}},\phi_{n_{k}})]$ 
converging in $\wQ$ to a critical point $[(A,\phi)]$ of $\sw$. Moreover, 
$$\lim_{n_{k}\rightarrow\infty}\sw([(A_{n_{k}},\phi_{n_{k}})])=\sw([(A,\phi)])$$

\end{proposition} 

Consequently, the basic Deformation Lemma of Morse Theory can be applied 
allowing the application of the \emph{Minimax Principle}(~\cite{Pa68}).

\section{\bf{W-Homotopy Type of $\wQ$}}

In this section, the hypothesis of the theorem ~\ref{Th:5018} are checked to the space $\Q$, and
the study of the weak homotopy type of $\Q$ is performed.

 We start observing that;
\begin{enumerate}
\item The quotiente spaces $\B =\Aa\diagup \G$ and $\vsa\diagup \G$ are Hausdorff spaces (~\cite{FU84}).
\item  the $\G$-action on $\Aa$ is not free since the action of the subgroup  of constant maps  
$g:M\rightarrow U_{1}$, $g(x)=g, \forall x\in M$, acts trivially on $\Aa$. 
\end{enumerate}

As mentioned before, instead of the $\G$-action, we consider on $\ca$ the  $\widehat{\G}$-action. 
On $\Aa$, the $\widehat{\G}$-action is free, and so, the space 
$\widehat{\mathcal{B}}_{\alpha}=\Aa/\wG$ is a 
manifold. 

The $\wG$-action on $\vsa$  is  free except on the 0-section, where the 
isotropic group is the full group $\wG$. The action also preserves the spheres in $\vsa$, 
consequently, the quotient space is a cone over the quotient of a sphere by the $\wG$-action.
 Therefore, the quotient space is contractible.

It follows from the Corollary of ~\ref{P:99} that  there exists the fibration

$$\vsa \rightarrow \Aa\times_{\wG}\vsa \rightarrow \widehat{\mathcal{B}}_{\alpha}$$

By the contractibility of $\vsa$, it follows that 

$$\Aa\times_{\wG}\vsa \overset{\text{htpy}}{\sim}\widehat{\mathcal{B}}_{\alpha} $$

In ~\cite{AJ78}, they studied the homotopy type of the space 
$\mathcal{A}\diagup\mathcal{G^{*}}$,  where $\mathcal{A}$ is the
space of connections defined on a G-Principal Bundle P and

$$\mathcal{\mathcal{G}^{*}}=\{g\in \mathcal{G} \mid g(x_{0})=I\},$$

\noindent is a subgroup og the Gauge Group  $\mathcal{G}=\varGamma(Ad(P))$. 
They observed that $\mathcal{G}^{*}$ acts freely on $\mathcal{A}$, and so, the
quotiente space $\mathcal{B}^{*}$ is a manifold. We need to compare the 
$\wG$ and $\G^{*}$ actions, neverthless, they turn out to be equal. The exact sequence

$$1\rightarrow U_{1}\rightarrow \G \overset{\rho}{\rightarrow}\G^{*}\rightarrow 1,
\quad \rho(g)=g(x_{0})^{-1}.g$$ 

\noindent implies that $\G^{*}\overset{\text{diffeo}}{\simeq} \G/U_{1}=\wG$,, where the quotient 
$U_{1}$ corresponds to the constant maps in $Map(X,U_{1})$. The actions are equal.

In this way, the  results of ~\cite{AJ78} can be applied to the understanding
of  the topology of the space $\Aa/\wG$.

 The weak homotopy type of $\B^{*}$ has been studied in ~\cite{AJ78} and ~\cite{DK91};
they proved the following; 

\begin{theorem}
Let $\mathcal{L}_{\alpha}$ be a complex line with $c_{1}(\mathcal{L}_{\alpha})=\alpha$, 
$\mathcal{E}U_{1}$ be the Universal bundle associated to $U_{1}$  and 

$$ Map^{0}_{\alpha}(X,\C P^{\infty})=\{f:X\rightarrow \C P^{\infty}\mid 
f^{*}(\mathcal{E}U_{1})\overset{\text{iso}}{\sim} \la , f(x_{0})=y_{0}\}.$$

 Then,

$$\B^{*} \overset{\text{w-htpy}}{\sim} Map^{0}_{\alpha}(X,\C P^{\infty})$$

\end{theorem}


\begin{corollary}
The space $\wQ$ is path-connected and

$$\pi_{n}(\wQ)=\pi_{n}(Map^{0}_{\alpha}(X,\C P^{\infty})),\quad n\in\N.$$
\end{corollary}

The set of path-connected components of $Map^{0}(X,\C P^{\infty})$ is equal to the
space  of homotopic classes $f:X\rightarrow\C P^{\infty}$, denoted by $[X,\C P^{\infty}]$. 
From Algebraic Topology, we know that 
\begin{enumerate}
\item There is a 1-1 correspondence
$$\{\mathcal{L}\mid \mathcal{L}\quad\text{is a complex line bundle over X}\}\leftrightarrow 
Map^{0}(X,\C P^{\infty}),$$
 \item The space of isomorphic classes of complex line bundles is 1-1 with
  $[X,\C P^{\infty}]$, i.e., if $\mathcal{L}$ is isomorphic to $\la$ then 
$f \in Map^{0}_{\alpha}(X,\C P^{\infty})$.
\item $[X,\C P^{\infty}]=H^{2}(X,\Z)$.
\end{enumerate}
 
In other words, 

$$\pi_{0}(Map^{0}(X,\C P^{\infty}))=H^{2}(X,\Z)$$

 
By the \emph{Minimax Principle}, the $\sw$-functional attains its minimum value 
in $\wQ$, and so, the equations (~\ref{E:09}) and (~\ref{E:10}) admit a solution. 

\begin{theorem}
Let $\alpha\in Spin^{c}(X)$. For each $n\in\N$, the homotopy group 
$\pi_{n}(Map^{0}(X,\C P^{\infty}))$, $n\in\N$
is isomorphic to 

$$\mathcal{H}= H^{2}(X,\Z)\oplus \{H^{1}(S^{n},\Z)\otimes H^{1}(X,\Z)\}
\oplus H^{2}(S^{n},\Z),$$

\noindent and it is computed in table 1.
\end{theorem}
\begin{proof}
Since 

$$\pi_{n}(Map^{0}(X,\C P^{\infty})) \simeq H^{2}(S^{n}\times X,\Z)$$

\noindent we can perform the computation of $\pi_{n}(Map^{0}_{\alpha}(X,\C P^{\infty})) $
by fixing a class of  $[X,\C P^{\infty}]$.

For a class $\alpha \in H^{2}(X,\Z)$, we fix a map $f:X\rightarrow \C P^{\infty}$ representing 
$\alpha$ and  $a\in S^{n}$. Thus, 

$$\pi_{n}(Map^{0}_{\alpha}(X,\C P^{\infty}))= 
[(S^{n}\times X,\{a\}\times X),(\C P^{\infty}, f(x_{0}))]$$

$$= [(S^{n}\times X,\{a\}\times X) , \C P^{\infty}]$$

However,
$$[(S^{n}\times X, \{a\}\times X) , \C P^{\infty}] = 
H^{2}(S^{n}\times X,\Z)/ H^{2}(\{a\}\times X,\Z)$$

Let $\mathcal{H}= H^{2}(S^{n}\times X,\Z)/ H^{2}(\{a\}\times X,\Z)$.

By Kuneth's formula,

$$H^{2}(S^{n}\times X,\Z)=H^{2}(X,\Z)\oplus \{H^{1}(S^{n},\Z)\otimes H^{1}(X,\Z)\}
\oplus H^{2}(S^{n},\Z)$$

\vspace{05pt}

\noindent and $H^{2}(\{a\}\times X ,\Z)=H^{2}(X,\Z)$. Consequently,

$$\mathcal{H}=\{H^{1}(S^{n},\Z)\otimes H^{1}(X,\Z)\}\oplus H^{2}(S^{n},\Z)$$

The group $\mathcal{H}$ is described below in the \emph{Table 1};
(the symbol $*$ stands whether the group is 0 or not)

\begin{table}[htbp]
  \centering
  \setlength{\belowcaptionskip}{10pt}
  \caption{$\mathcal{H}$ }
  \label{tab:hor-extern}
  \begin{tabular}{|l|l|l|}
    \hline
    \hline
    \emph{n}   &    {$H^{1}(X,\Z)$}         & $\mathcal{H}$  \\
    \hline
    $\ne$ 1, 2    &   *          &         0 \\
    1             &        0       &         0 \\
    1             &      $\ne 0$   &     $H^{1}(X,\Z)$ \\
    2             &          *   &     $\Z$ \\
    \hline
    \hline
  \end{tabular}
\end{table}
\end{proof}

\begin{theorem} \label{Th:M}
There exist critical points for the $\sw$-functional. 
\end{theorem}
\begin{proof}
The \emph{Minimax Principle} implies that the minimum value is attained in all 
connected component of $\Q$. The table 1 shows that it is possible to construct 
non-contractible families of elements in $\Q$, consequently, by applying the 
\emph{Minimax Principle}, it follows that there exist stable and non-stable 
critical points. In other words, there exist  solutions to the equations 
(~\ref{E:09}) and (~\ref{E:10}).
\end{proof}

\appendix

\section{Quotient Spaces by the Diagonal Action}

Let  M, N be smooth manifolds endowed with  G-actions
$\alpha_{M},\alpha_{N}$ (respec.). About the G-actions, we will assume that;
\begin{enumerate}
\item The isotropic groups of the action on M are isomorphics, i.e. there 
exists a Lie Group H such that for all $m\in M$ $\Rightarrow$ $G_{m}\simeq H$, 
\item The quotient spaces $M\diagup G$ and $N\diagup G$ are Hausdorff spaces.  
\end{enumerate}

The product action of $G\times G$ on the manifold $M\times N$,  is defined  by   

$$\alpha_{M}\times\alpha_{N}:G\times G\times(M\times N)\rightarrow M\times N,$$

$$\alpha_{M}\times\alpha_{N}(g_{1},g_{2},m,n)=(\alpha_{M}(g_{1},m),
\alpha_{N}(g_{2},n)),$$ 

\noindent or equivalently,

$$(g_{1},g_{2}).(m,n)=(g_{1}.m,g_{2}.n)$$

\begin{definition}
The diagonal action $\alpha_{\mathcal{D}}:G\times (M\times N)\rightarrow
M\times N$ is defined as 
$$\alpha_{\mathcal{D}}(g,(m,n))=(\alpha_{M}(g.m),\alpha_{N}(g,n)),$$
\noindent and denoted as g.(m,n)=(g.m,g.n). The quotiente space is denoted by
$M\times_{G} N$.

\end{definition}

\begin{definition}
Let $m\in M$ and $n\in N$. The corresponding  orbit  spaces are defined as follows;
\begin{enumerate}
\item For the action $\alpha_{M}$ on M, let 
$\mathcal{O}^{M}_{m}=\{g.m\mid g\in G\}$;
\item For the action $\alpha_{N}$ on N, let $\mathcal{O}^{N}_{n}=\{g.n\mid
  g\in G\}$;
\item For the product action ($P$-action) $\alpha_{M}\times\alpha_{N}$ on $M\times N$, let
  $$\mathcal{O}^{P}_{(m,n)}=\{(g_{1}.m,g_{2}.n)\mid g_{1}, g_{2}\in G \};$$
\item For the diagonal action $\mathcal{D}$-action)  $\alpha_{\mathcal{D}}$ on $M\times N$, let 
$$\mathcal{O}^{\mathcal{D}}_{(m,n)}=\{(g.m,g.n)\mid g\in G\}$$
\end{enumerate}
\end{definition}

The orbit of $(m,n)$, by the product action, is easily described by the orbits in M and N as 
 
$$\mathcal{O}^{P}_{(m,n)}=\mathcal{O}^{M}_{m}\times \mathcal{O}^{N}_{n}.$$

\noindent Consequently,

$$(M\times N)\diagup (G\times G)=(M\diagup G)\times (N\diagup G),$$

\noindent what induces the fibration

$$N\diagup G \longrightarrow (M\times N)\diagup (G\times G) \longrightarrow M\diagup G$$

\noindent


In order to describe the topology of the space $M\times_{G} N$, we 
 consider the commuting diagram

\begin{equation}
\begin{CD}
M\times N @>p_{1}>> M \\
@V\text{$\pi^{M\times N}$}VV @V\text{$\pi^{M}$}VV\\
M\times_{G} N @>\mathfrak{p}>> M\diagup G
\end{CD}
\end{equation}

\vspace{10pt}

\noindent where 
\begin{enumerate}
\item $p_{1}:M\times N\rightarrow M$ is the projection on the 1$^{st}$ factor;
\item $\pi^{M\times N}:M\times N\rightarrow  M\times_{G} N$ is the
  projection induced by the quotient;
\item $\mathfrak{p}: M\times_{G} N \rightarrow M\diagup G$ is the 
natural map induced by the projection
$\mathcal{O}^{\mathcal{D}}_{(m,n)}\rightarrow \mathcal{O}^{M}_{m}$.
\end{enumerate}
 
 From now on, we fix $[m_{0}]\in M\diagup G$ in order to describe 
$\mathfrak{p}^{-1}([m_{0}])$.   

From the diagram, we get that
\begin{enumerate}
\item $(\pi^{M})^{-1}([m_{0}])=\mathcal{O}^{M}_{m_{0}}$
\item $(\pi^{M}\circ p_{1})^{-1}([m_{0}]) = \mathcal{O}^{M}_{m_{0}}\times N$,
\item $(\pi^{M\times N})^{-1}([m_{0},n_{0}])=
\mathcal{O}^{\mathcal{D}}_{[(m_{0},n_{0})]}$
\end{enumerate}

\vspace{15pt}

\begin{proposition}
The subspace $\mathcal{O}^{M}_{m_{0}}\times N$ is a G-space with respect to the $\mathcal{D}$-action.
\end{proposition}

\begin{proof}
The proof is splited in two easy claims;
\begin{enumerate}
\item If $(m,n)\in \mathcal{O}^{M}_{m_{0}}\times N$, then 
$\mathcal{O}^{\mathcal{D}}_{(m,n)}\subset \mathcal{O}^{M}_{m_{0}}\times N$.

Let $m=g.m_{0}$;

$$g^{,}.(m,n)= (g^{,}g.m_{0},g^{,}.n)\in \mathcal{O}^{M}_{m_{0}}\times N$$

\item If $(m,n)\in \mathcal{O}^{M}_{m_{0}}\times N$, then there exits $g\in G$ 
and $n^{,}\in N$ such that $(m,n)\in \mathcal{O}^{\mathcal{D}}_{(m_{0},n^{,})}$

Let $m=g.m_{0}$ and $n^{,}=g^{-1}.n$;

$$(m,n)=g.(m_{0},g^{-1}.n)\quad \Rightarrow (m,n)\in\mathcal{O}^{\mathcal{D}}_{(m_{0},n^{,})}$$

\end{enumerate}
\end{proof}

Consequently, 

$$\mathfrak{p}^{_1}([m_{0}])= \mathcal{O}^{M}_{m_{0}}\times_{G} N$$ 

\begin{proposition}\label{P:01}
$$\mathcal{O}^{\mathcal{D}}_{(m,n)}\cap p^{-1}_{1}(m_{0})=\{g.n\mid g\in G_{m_{0}}\}$$
\end{proposition}
\begin{proof}
Let $m=g.m_{0}$;  so $(m,n)=g.(m_{0},g^{-1}.n)$ $\Rightarrow$ 
$\mathcal{O}^{\mathcal{D}}_{(m,n)}= \mathcal{O}^{\mathcal{D}}_{(m_{0},g^{-1}.n)}$. 

Neverthless,
 
$$g.(m_{0},n)\in p^{-1}_{1}(m_{0}) \Leftrightarrow \exists g\in G\quad \text{such that}\quad 
g.(m,n)=(m_{0},n^{,}),$$

\noindent this implies that $g\in G_{m_{0}}$ and $n^{,}=g.n$
\end{proof}

Therefore, every  $\mathcal{D}$-orbit meet the set $p^{-1}_{1}(m_{0})$. Now, we will construct a 
smooth invariant map. Define  $\rho:\mathcal{O}^{M}_{m_{0}}\times N\rightarrow N/G_{m_{0}}$  by 

$$\rho(g^{,}.(g.m_{0},n))=g^{-1}.n,\quad\text{for all}\quad g^{,}\in G$$

The map $\rho$ is well defined, since whenever 
$g\in G_{m_{0}}$ it follows that

$$g^{-1}.n = \rho((g.m_{0},n)) = \rho((m_{0},n)) = n $$

This remark is consistent with ~\ref{P:01}.

It is easily seen to be bijective. Consequently,

\begin{proposition}\label{P:99}
$$M\times_{G} N=\bigcup_{[m]\in M/G}N/G_{m}$$
\end{proposition}
\begin{proof}
Its follows from the discussion above, since it was concluded that

$$\mathfrak{p}^{-1}([m_{0}])= N/G_{m_{0}}$$

\end{proof}

\begin{corollary}
\begin{enumerate}
\item If the G-action on M is free, then there is the fibration

\begin{equation}\label{F:01}
N \longrightarrow M\times_{G} N \longrightarrow M/G
\end{equation}

\item Suppose that there exist a Lie Group $H$ such that for all $m\in M$ it is true that
$G_{m}$ is isomorphic to $H$. In this case, there is also a fibration

\begin{equation}
N/H \longrightarrow M\times_{G} N \longrightarrow M/G
\end{equation}

\noindent which fiber $N/H$  may be is a singular space.
\end{enumerate}
\end{corollary}

\vspace{20pt}

\noindent\emph{Universidade Federal de Santa Catarina \\
    Campus Universitario , Trindade\\
               Florianopolis - SC , Brasil\\
               CEP: 88.040-900}
\vspace{10pt}

\end{document}